\documentclass[12pt,a4paper,utf8]{article}
\usepackage[margin=1in]{geometry}  
\usepackage{graphicx}              
\usepackage{amsmath}               
\usepackage{amsfonts}              
\usepackage{amsthm}                
\usepackage{amsfonts,amsmath,amssymb,euscript,amsthm,indentfirst,latexsym,delarray,hhline,graphicx}
\usepackage{hyperref}
\usepackage{epstopdf}
\usepackage{hyperref}
\usepackage{cmap}
\usepackage{mathtext}
\usepackage{type1ec}
\newtheorem*{remark}{Remark}
\newtheorem{lemma}{Lemma}

\newtheorem*{example}{Example}
\newtheorem{proposition}{Proposition}

\begin{document}

\title{ Construction of an algebra corresponding to a statistical model of the square ladder (square lattice with two lines) }

\author{Valerii Sopin}

\maketitle

\begin{abstract}
In this paper we define infinite-dimensional algebra and its
representation, whose basis is naturally identified with
 semi-infinite configurations of the square ladder model.

We also extrapolate the ideas for the cyclic 3-leg triangular
ladder. All of these propose a way for generalization, which leads
to representations of $N=2, \dots$ algebras.

\end{abstract}

\section{Introduction}

For each graph $\Gamma$ we can construct a statistical model in
which the set of configurations is the set of arrangements of
particles at graph vertices such that at each vertex at most one
particle is located and two particles cannot be located at vertices
joined by an edge.

The previous paper $[1]$ discussed combinatorial properties of the
set of configurations of the $2\times n$ square lattice (or simply
the square ladder) graph:
\begin{figure}[h]
\centering
\includegraphics[width=10.0cm]{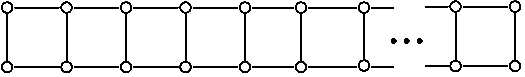}
\caption{$2\times n$ square lattice model} \label{fig.0}
\end{figure}

Let's assign the fermion algebra of anti-commuting elements $x_i$
and $y_i$ to the graph in Fig. 1, where  $i=-k, \dots, k$ if $n =
2k+1$ or $i=-k, \dots, k-1$ if $n = 2k$; moreover, elements $x_i,
y_i$ satisfy the next relations: $$x_iy_i=0,\; x_ix_{i+1}=0,\;
y_iy_{i+1}=0.$$

Studies of integrable models of statistical mechanics show a close
connection between the set of configurations of the corresponding
lattice graph and the representation of some infinite-dimensional
algebra $[2][3][4][5][6][7][8]$. Based on this, the paper yields the
following results:

In Section 2, based on idea of semi-infinite forms, the set of
configurations is defined for the square ladder model, which is
infinite in both directions. A bigraduation is introduced and
statistical sum is calculated.

In Section 3, such a deformation of the fermion algebra for the
graph of the square ladder model, which is infinite in both
directions, is determined (the obtained algebra is close to
conformal algebras) that the character of its representation is
equal to the statistical sum from Section 2.

In Section 4, the cohomology of complexes, constructed from
finite-dimensional quotient algebras of the deformed algebra from
Section 3, are calculated. The corresponding complexes are either
acyclic, or their cohomology is one-dimensional. The deformation was
selected by the latter, among other things.

In Section 5, further generalization is discussed.

\section{Set of semi-infinite configurations}

Instead of the square ladder model, Fig. 1, it is convenient to
consider an equivalent combinatorial problem (the results of paper
$[1]$ remain valid), the twisted square ladder model:

\begin{figure}[h]
\centering
\includegraphics[width=10cm]{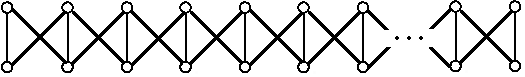}
\caption{twisted square ladder model} \label{fig.0}
\end{figure}

According to $[8]$, the set $\Omega \subset \mathbb{Z}$ is called a
Dirac set if $\Omega_e = \Omega \cap \mathbb{Z}_{\geq 0}$ and
$\Omega_p = \mathbb{Z}_{< 0} \setminus \Omega$ are finite. The value
$C(\Omega) = \text{count}(\Omega_e) - \text{count}(\Omega_p)$ is
called the charge of the set $\Omega$, and $U(\Omega) =
\sum\limits_{\alpha \in \Omega_e \cup \Omega_p} |\alpha|$ is called
its energy.

The set of such configurations $\{ \dots, v_{i_{n-1}},
v_{i_{n}},v_{i_{n+1}}, \dots, v_{i_{1}}\}$ of the graph of the
twisted square ladder model, which is infinite in both directions,
that the set of indexes $\Omega =\{i_n\}$ is a Dirac set, is called
the set of semi-infinite configurations $\Delta^{\frac{\infty}{2}}$.
The location of particle $v_{i_{n}}$ is given by $x_{i_{n}}$ or
$y_{i_{n}}$ depending on whether the particle is in the upper or
lower row, where $i_{n}$ is the corresponding column number.

It follows from the definition of Dirac set that a configuration has
``tail'' (see Fig. 3), i.e. there is such $N \in \mathbb{Z}$ that
$$i_{N-j+1}=i_{N}-j+1 \; \forall j \in \mathbb{N}.$$

Moreover, it follows from the definitions of set of configurations,
Dirac set and the type of the graph under consideration (Fig. 2)
that all $v_{i_{N-j+1}}$ are elements of one row.

\begin{figure}[h]
\centering
\includegraphics[width=9cm]{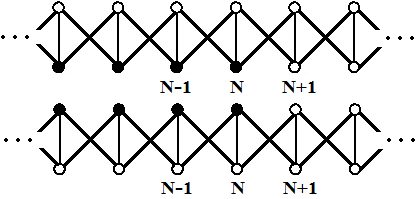}
\caption{tails of twisted square ladder model} \label{fig.0}
\end{figure}

We will talk about energy and charge of $\{ \dots, v_{i_{n-1}},
v_{i_{n}},v_{i_{n+1}}, \dots, v_{i_{1}}\}$ implying energy and
charge of the corresponding set of indexes $\Omega =\{i_n\}$.

\subsection{Calculation of statistical sum}

In this subsection, our goal will be to calculate the statistical
sum of the subset of configurations $\Delta^{\frac{\infty}{2}}$ with
the fixed type of tail (Fig. 3) relative to $C$ and $U$, i.e.
$$F(t, q) = \sum_{ across\; all\; Dirac\; sets\; \Omega }
t^{C(\Omega)}q^{U(\Omega)}.$$

\begin{proposition} $$F(t, q)= \prod\limits_{m=1}^{\infty} \frac{1+q^m}{1-q^m}\sum\limits_{n=-\infty}^{\infty}
 q^{n(n-1)/2} t^{n}.$$ \end{proposition}
\textit{Proof. }Let's consider
$$F_{+}(t, q) = \sum_{\substack{ across\; all\; Dirac\; sets\; \Omega \\ with\; \Omega_p = \varnothing }}
t^{C(\Omega)}q^{U(\Omega)} = \sum_{i=0}^{\infty} \varphi_n(q)t^i.$$

\begin{lemma} $$F_{+}(t, q)=\sum_{n=0}^{\infty} \frac{(q)^{+}_n}{(q)_n}
q^{n(n-1)/2} t^n,$$ where $(q)^{+}_n =(1+q)(1+q^2)\dots(1+q^n)$,
$(q)_n =(1-q)(1-q^2)\dots(1-q^n)$. \end{lemma} \textit{Proof. }One
can readily see that the following is true
$$F_{+}(t, q) = tF_{+}(qt, q)+\dot{F}_{+}(qt, q),$$ where
$$\dot{F}_{+}(qt, q)=2tqF_{+}(q^2t, q)+\dot{F}_{+}(q^2t, q),$$ as we consider all possible options
for particle arrangement at the column no. 0.

Then,
$$F_{+}(t, q) - tF_{+}(qt, q)= F_{+}(qt, q) - tqF_{+}(q^2t, q) + 2tqF_{+}(q^2t, q).$$

Therefore, we obtain functional equation:
$$F_{+}(t, q) = F_{+}(qt, q) + t(F_{+}(qt, q) + qF_{+}(q^2t, q)).$$

Hence,
$$\varphi_n(q) =
     q^n\varphi_n(q)+q^{n-1}\varphi_{n-1}(q)+q^{2n-1}\varphi_{n-1}(q).$$

Consequently,
$$(1-q^n)\varphi_n(q) = q^{n-1}(1+q^n)\varphi_{n-1}(q).$$

As $\varphi_0(q)=1$, it has been demonstrated$._\bigtriangledown$

Let's consider another statistical sum
$$F_{N}(t, q) = \sum_{\substack{\text{across all Dirac sets $\Omega$:} \\ \text{$\Omega_p$ contains only the numbers not exceeding -N}}}
t^{C(\Omega)}q^{U(\Omega)}.$$

Let's consider the transformation associating the Dirac set $\Omega
= \{i_k\}$ to the Dirac set $\Omega_N = \{i_k+N\}$. It is easily
seen that  $C(\Omega_1) = C(\Omega)+1$ and $U(\Omega_1) =
U(\Omega)+C(\Omega)+1$. Following this, by induction,
$$C(\Omega_N) = C(\Omega)+N \text{ and } U(\Omega_N) =
U(\Omega)+NC(\Omega)+N(N+1)/2.$$

A Dirac set contains all numbers not exceeding $-N$ if and only if
$(\Omega_{N-1})_p =  \varnothing$. This yields the equation
$$F_{+}(t, q) = t^{N-1}q^{N(N-1)/2}F_{N}(tq^{N-1}, q).$$

Consequently,
$$F_{N}(t, q) = \frac{q^{(N-1)(N-2)/2}}{t^{N-1}}F_{+}(t/q^{N-1}, q) =$$
$$ = \sum_{n=0}^{\infty} \frac{(q)^{+}_n}{(q)_n}
q^{n(n-1)/2+(N-1)(N-2)/2-n(N-1)} t^{n-N+1} = \sum_{n=-N+1}^{\infty}
\frac{(q)^{+}_{n+N-1}}{(q)_{n+N-1}} q^{n(n-1)/2} t^{n}.$$

Going to the limit $N\longrightarrow \infty$, we obtain the
required$._\bigtriangledown$

\begin{remark}$$\prod\limits_{m=1}^{\infty} \frac{1+q^m}{1-q^m} =
\prod\limits_{m=1}^{\infty} \frac{1}{1-q^m}
\prod\limits_{k=0}^{\infty} \frac{1}{1-q^{2k+1}} =$$
$$ = (\prod\limits_{m=1}^{\infty} \frac{1}{1-q^m})^2\prod\limits_{k=1}^{\infty}
(1-(q^2)^k),$$ due to the equality of generating functions for
splitting into odd and different summands $[9]$.
\end{remark}

\section{$\mathbf{SqL}$ algebra}

The fermion algebra for the infinite in both directions graph as per
Fig. 2 is the following algebra of anti-commuting elements $x_i$ and
$y_i$, $i \in \mathbb{Z}$:
$$\mathbb{C}[\dots, x_{-1}, y_{-1}, x_{0}, y_{0}, x_{1}, y_{1},\dots] / (x_iy_i=0,\; x_iy_{i+1}=0,\; x_{i+1}y_{i}=0).$$

Let's determine the deformation of this algebra, namely the
$\mathbf{SqL}$ algebra, generated by anti-commuting elements $x_i$,
$y_i$, $i \in \mathbb{Z}$, satisfying the relations below (let's
denote it as $\verb"REL"$):
$$\forall s \in \mathbb{Z} \sum_{i\in  \mathbb{Z}} x_i y_{-i+s} = 0,\; \sum_{i\in  \mathbb{Z}} x_i y_{-i+2s+1}(-1)^i= 0,$$
with additional action of two operators $c$ and $u$:
$$[c,\;u]=0; \;$$ $$[c,\;x_i]=x_i,\;[c,\;y_i]=y_i;\; $$ $$[u,\;x_i]=-ix_i,\;[u,\;y_i]=-iy_i.$$

Let's denote generating functions  $X(t) = \sum_{i\in  \mathbb{Z}}
x_i t^{-i}$, $Y(t) = \sum_{i\in \mathbb{Z}} y_i t^{-i}$, then the
relations $\verb"REL"$ can be rewritten as

$$X(t)Y(t) = 0,\; X(t)Y(-t)= X(-t)Y(t).$$

The relations $\verb"REL"$ are infinite; therefore, formally, the
algebra with such relations has no sense. However, if we consider
its representation in a graduated space, with peak limiting for the
graduation, as is customary in the theory of conformal algebras,
then everything will be determined.

Accordingly, for any integer $N$ let's define induced
representations with extreme vectors, i.e. spaces $\Upsilon(N)$
spawned by elements $x_{i}, y_{i},$ $i \in \mathbb{Z}$, from the
 vectors $\gamma_N$, for which the following is true:
$$x_{i\geq N-1}\circ\gamma_{N}=0,\; y_{i\geq N}\circ\gamma_{N}=0.$$

In other words, $\Upsilon(N) = (\mathbb{C}[\dots, x_{-1}, y_{-1},
x_{0}, y_{0}, x_{1}, y_{1},\dots] / \verb"REL") \circ\gamma_{N}$.

There are mappings
$$\dots \longrightarrow \Upsilon(-2) \longrightarrow \Upsilon(-1)\longrightarrow \Upsilon(0)\longrightarrow \Upsilon(1) \longrightarrow \dots,$$
which follow from $\gamma_{N-1} \longmapsto y_{N}\circ\gamma_{N}.$

Mappings $\Upsilon(N) \longrightarrow \Upsilon(N+1)$ are embeddings.
This result is derived from the existence of monomial basis (Lemma
2) in $\Upsilon(N)$. Now we are going to discuss it. But, to begin
with, let's construct an auxiliary representation of the
$\mathbf{SqL}$ algebra by using the Clifford algebra.

\subsection{Constructing representation of $\mathbf{SqL}$ algebra}

Let \begin{equation*} \delta_{i+j} = \begin{cases} 1, &i+j=0;\\ 0,
&\text{otherwise.} \end{cases} \end{equation*}

The algebra  $OCl$ is generated by elements $o_{-2k-1}, k \in
\mathbb{N} \cup \{0\},$ $\varphi_{i}, \varphi^{*}_{i}, i \in
\mathbb{Z}$, for which the following is true:
$$[o_{-2k-1}, o_{-2l-1}] =0, [o_{-2k-1}, \varphi_{i}] =0, [o_{-2k-1}, \varphi^{*}_{i}] =0,$$
$$[\varphi_{i}, \varphi_{j}]_{+}=0, [\varphi^{*}_{i},
\varphi^{*}_{j}]_{+}=0,[ \varphi_{i},\varphi^{*}_{j}]_{+}=\delta_{i+j}.$$

Moreover, additional action of $s$ and $w$ is defined:
$$[s,\;w]=0, \;$$
$$[s,\;\varphi_i]=\varphi_i,\;[s,\;\varphi^{*}_i]=-\varphi^{*}_i,\;[s,\;o_{-2i-1}]=0,$$
$$[w,\;\varphi_i]=-i\varphi_i,\;[w,\;\varphi^{*}_i]=-i\varphi^{*}_i,\;[w,\;o_{-2i-1}]= (2i+1) o_{-2i-1}.$$

Let's introduce generating functions: $$\varphi(z) =
\sum\limits_{i\in \mathbb{Z}} \varphi_{i} z^{-i}, o(z)
=\sum\limits_{i\in \mathbb{N} \cup \{0\}} o_{-2i-1} z^{2i+1},
\varphi^{*}(z) = \sum\limits_{i\in \mathbb{Z}} \varphi^{*}_{i}
z^{-i}.$$

Let's define: $$\psi(z)= o(z) \varphi(z) = \sum\limits_{i\in
\mathbb{Z}} \psi_{i} z^{-i}, \text{ where } \psi_{i}
=\sum\limits_{j+l=i}o_{j}\varphi_l.$$

The relations $\verb"REL"$ are fulfilled for $\varphi(z)$ and
$\psi(z)$:
$$\varphi(z)\psi(z) = o(z) \varphi^2(z) = 0,$$ $$\varphi(z)\psi(-z) =
\varphi(z)\varphi(-z)o(-z)= -\varphi(z)\varphi(-z)o(z)
=\varphi(-z)\psi(z).$$

Let's suppose there exists extreme vector $v_0$:
$$\varphi_{i>0} \circ v_0 =0, \varphi^{*}_{i\geq0} \circ v_0 =0.$$

Let's set an action for $s$ and $w$ on $v_0$:  $$s\circ v_0=0,\; w
\circ v_0=0.$$

The character  of $ OCl \circ v_0 = \mathbb{C}[\varphi_{0},
\varphi^{*}_{-1},o_{-1}, \varphi_{-1}, \dots] \circ v_0$ is defined
as $Tr_{OCl \circ v_0}(q^{w}t^{s}).$

\begin{proposition} $Tr_{OCl \circ v_0}(q^{w}t^{s})$ and $F(t, q)$ coincide. \end{proposition}
\textit{Proof. }As is well known, the character of the space of
polynomials in anti-commuting variables $\varphi_{i},
\varphi^{*}_{i-1}, i=0, -1, \dots,$ is equal to
$$\prod\limits_{j=0}^{\infty} (1+q^{j}t) \prod\limits_{k=1}^{\infty} (1+q^{k}t^{-1}),$$ which equals to $$(\prod\limits_{l=1}^{\infty} \frac{1}{1-q^l}) \sum\limits_{n=-\infty}^{\infty}
 q^{n(n-1)/2} t^{n}$$ according to the Jacobi triple product $[9]$.

Moreover, as is well known, the character of the space of
polynomials in commuting variables $o_{2i-1}, i=0, -1, \dots,$ is
equal to $$\prod\limits_{m=0}^{\infty}
\frac{1}{1-q^{2m+1}}._\bigtriangledown$$

\subsection{Direct limit}

\begin{lemma} In $\Upsilon(N)$ there is a monomial basis consisting
of $$x_{i_1}x_{i_2}\dots x_{i_l}y_{j_1}y_{j_2}\dots
y_{j_m}\circ\gamma_{N}, |i_{a}-j_{b}|\geq2.$$\end{lemma}
\textit{Proof. }Using $\verb"REL"$, as well as skew symmetry, any
monomial $x_{\widetilde{z}_1}\dots
x_{\widetilde{z}_k}y_{\widetilde{t}_1}\dots y_{\widetilde{t}_l}$ can
be expressed in the required form as above (see $[7]$ and $[8]$):
each of newly appeared monomials will be lexicographically less than
the original ($\text{deg}(x_j) = \text{deg}(y_j) = j$, all $y_j$
follow after all $x_i$). Of course, newly appeared monomials may
contain pairs breaking the conditions $|i_{a}-j_{b}|\geq2$, but the
lexicographical order allows to develop an iterative procedure.

Let's demonstrate the linear independence. For this purpose, let's
use $\mathbf{SqL}$ representation defined above. Let's set the
homomorphism of algebras: $$x_i \longmapsto \varphi_{i},\; y_i
\longmapsto \psi_{i},\;u\longmapsto w,\; c \longmapsto s,$$

from here we obtain the mapping of spaces:
$$\Upsilon(N) \longrightarrow \Upsilon_{\varphi, \psi}(N)=\mathbb{C}[\dots, \varphi_{-1}, \psi_{-1}, \varphi_{0}, \psi_{0}, \varphi_{1}, \psi_{1},\dots] \circ \gamma_N.$$

By definition,  $\psi_{i} =\sum\limits_{j+l=i}o_{j}\varphi_l$ acts
with non-zero to $\gamma_N$ only with a finite number of
$o_{j}\varphi_l,$  therefore:
$$\Upsilon_{\varphi, \psi}(N) = \mathbb{C}[\dots, \varphi_{-2}, o_{-2},\varphi_{-1}, o_{-1}, \varphi_{0}, \varphi_{1}, \varphi_{2},\dots] \circ \gamma_N.$$

By definition of  $\gamma_N$, the latter is
$$\mathbb{C}[\dots, \varphi_{-2}, o_{-2},\varphi_{-1}, o_{-1}, \varphi_{0}][\varphi_{1}, \varphi_{2},\dots,  \varphi_{N-1}] \circ \gamma_N.$$

Then, from the definition of $OCl$, we obtain the following
isomorphism:
\begin{equation*}
\Upsilon_{\varphi, \psi}(N)\simeq
 \begin{cases}
     \mathbb{C}[\varphi_{0}, o_{-1}, \varphi_{-1}, o_{-3}, \dots][\varphi^{*}_{-1}, \varphi^{*}_{-2},\dots,  \varphi^{*}_{-N+1}] \circ v_0, &\text{$N >0$,}\\
   \mathbb{C}[\varphi_{0}, o_{-1}, \varphi_{-1}, o_{-3}, \dots] \circ \varphi_{N}\dots\varphi_{0} \circ v_0, &\text{$N \leq 0$.}
 \end{cases}
\end{equation*}

Therefore,
$$OCl \circ v_0 = \mathbb{C}[\varphi_{0}, \varphi^{*}_{-1},o_{-1}, \varphi_{-1}, \dots] \circ v_0 = \varinjlim \Upsilon_{\varphi, \psi}(N).$$

Thereafter, there are no additional relationships for
 $x_{i_1}x_{i_2}\dots x_{i_l}y_{j_1}y_{j_2}\dots
y_{j_m}\circ\gamma_{N},$ where $|i_{a}-j_{b}|\geq2$, because
otherwise the character $Tr_{OCl \circ v_0}(q^{w}t^{s})$ would be
less than the statistical sum $F(t, q)$, but from Proposition 2 we
know that this is not so.$_\bigtriangledown$

Let's define $$\Upsilon = \varinjlim \Upsilon(N).$$

Any vector of space $\Upsilon$ is a finite linear combination of
expressions
$$x_{i_1}x_{i_2}\dots x_{i_l}y_{j_1}y_{j_2}\dots y_{j_m}\circ\gamma_{N} = x_{i_1}x_{i_2}\dots x_{i_l}y_{j_1}y_{j_2}\dots y_{j_m}y_{N}\circ\gamma_{N+1}=\dots,$$
$$|i_{a}-j_{b}|\geq2,$$
for a sufficiently large $N$. Let's send $N$ to infinity, i.e.
formally substitute $\gamma_{N}$ with expression
$$y_{N}y_{N+1}y_{N+2}\dots\circ\gamma_{\infty}.$$

So, we obtain the representation, which can be naturally identified
with the set of semi-infinite configurations with the fixed type of
tail.

Let's set an action of  $c$ and $u$ on extreme vectors $\gamma_N$:
$$c \circ \gamma_N= -N\gamma_N,$$
\begin{equation*}
u \circ \gamma_N=
 \begin{cases}
     \frac{N(N+1)}{2}\gamma_N, &\text{$N \geq 0$,}\\
   \frac{N(N-1)}{2}\gamma_N, &\text{$N < 0$.}
 \end{cases}
\end{equation*}

The character of $\Upsilon$ is defined as
$Tr_{\Upsilon}(q^{u}t^{c}).$

\begin{proposition}  $Tr_{\Upsilon}(q^{u}t^{c})$ and $F(t, q)$ coincide. \end{proposition}
\textit{Proof. }After all the above it's obvious.$_\bigtriangledown$

\section{Cohomology}

Let's define algebra $A_n$ of anti-commuting elements $x_i$ and
$y_i$, where $i=-k, \dots, k$ if $n = 2k+1$ or $i=-k, \dots, k-1$ if
$n = 2k$, satisfying the relations $\verb"REL"$ (all other $x_j$ and
$y_j$, with other indexes, are assumed equal to zero).

Therefore, algebra $A_n$ is a finite-dimensional quotient algebra of
$\mathbf{SqL}$. It follows from Lemma 2 that the dimension of $A_n$
is identical to the dimension of the fermion algebra of the twisted
square ladder with $n$ columns. Moreover, they have the same
monomial basis.

Symbols $x_i$ or $y_i$ also denote the operator of multiplication by
the corresponding element.

Let's construct a complex  $K_n\cong A_n$ with differential
$x_0+y_0$.

Let's introduce elements from  $K_{2k+1}$, which are defined
recursively as follows:

\begin{minipage}[h]{.45\linewidth}
\begin{equation*}
\underline{h}_{2k+1}=
 \begin{cases}
     x_{-k} \underline{h}_{2k-3} y_{k}, &\text{$k \mod 2 \neq 0$,}\\
   y_{-k} \underline{h}_{2k-3} y_{k}, &\text{$k \mod 2 = 0$,}
 \end{cases}
\end{equation*}
\end{minipage}
\hfill
\begin{minipage}[h]{.45\linewidth}
\begin{equation*}
\overline{h}_{2k+1}=
 \begin{cases}
     y_{-k} \underline{h}_{2k-3} x_{k}, &\text{$k \mod 2 \neq 0$,}\\
   x_{-k} \underline{h}_{2k-3} x_{k}, &\text{$k \mod 2 = 0$,}
 \end{cases}
\end{equation*}
\end{minipage}\\
where $\underline{h}_3 = x_{-1}y_1, \overline{h}_3 = y_{-1}x_1;$
$\underline{h}_1 = y_0, \overline{h}_1 =y_0,$ $k \in \mathbb{N}
\bigcup \{0\}$.

\begin{figure}[h]
\centering
\includegraphics[width=16cm]{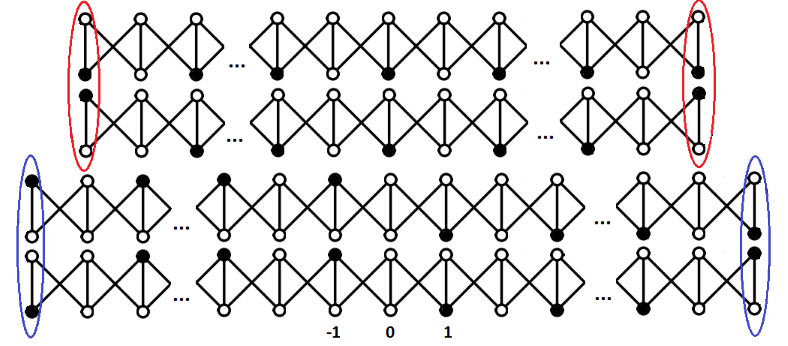}
\caption{$\underline{h}_{2k+1}$ and $\overline{h}_{2k+1}$}
\label{fig.0}
\end{figure}

\begin{proposition} The cohomology of the complexes $ K_{2k+1},\;K_{2k+2} $ is
one-dimensional. The element $\underline{h}_{2k+1}$ is a
representative of the corresponding cohomoloogy class in both
cases.\end{proposition} \textit{Proof. }To begin with, let's
consider the odd case, $ K_{2k+1}$.

Note that if $2s=\pm2k$, from $\sum_{i\in  \mathbb{Z}} x_i y_{-i+2s}
= 0$ it follows that $$x_{-k}y_{-k}=0=x_ky_k.$$

In addition, if $2s+1=-2k+1, 2k-1$, it follows from $\sum_{i\in
\mathbb{Z}} x_i y_{-i+2s+1} = 0$ and $\sum_{i\in \mathbb{Z}} x_i
y_{-i+2s+1}(-1)^i= 0$ that
 $$x_{-k}y_{-k+1}=x_{-k+1}y_{-k}= 0,\; x_{k}y_{k-1}=x_{k-1}y_{k} =0.$$

Then, let's introduce auxiliary complexes denoted as $K^{\alpha,
\beta}_{2k+1}, \alpha, \beta \in \{ x, y, \square \}.$

$K^{x, x}_{2k+1}$ means the image of operator of multiplication by
$x_{-k}x_{k}:$ $x_{-k}x_{k}K_{2k+1},$ $K^{x, y}_{2k+1}$ means the
image of $x_{-k}y_{k}:$ $x_{-k}y_{k}K_{2k+1}.$ $K^{x,
\square}_{2k+1}$ means the factor space of the image of operator of
multiplication by $x_{-k}:$ $x_{-k}K_{2k+1}/\{x_{k}=0,\;y_{k}=0\},$
$K^{\square, x}_{2k+1}$ is $x_{k}K_{2k+1}/\{x_{-k}=0,\;y_{-k}=0\}.$

Similar for $K^{y, y}_{2k+1}$, $K^{y, x}_{2k+1}$, $K^{y,
\square}_{2k+1}$, $K^{\square, y}_{2k+1}$.

$K^{\square, \square}_{2k+1}$ means
$K_{2k+1}/\{x_{-k}=0=y_{-k},\;x_{k}=0=y_{k}\} = K_{2k-1}.$

Therefore, we split $K_{2k+1}$, see Lemma 2 and $[1]$.

\begin{lemma} The complexes $K^{x, \square}_{2k+1}, K^{y,
\square}_{2k+1}, K^{\square, x}_{2k+1},K^{\square, y}_{2k+1}$ are
acyclic.

The cohomology of the complexes $\text{K}^{x, x}_{2k+1}$,
$\text{K}^{y, y}_{2k+1}$ is one-dimensional in case of odd $k$, and
elements $\overline{h}_{2k+1}$, $\underline{h}_{2k+1}$  are a
representatives of the corresponding cohomology class, respectively.
If $k$ is even, then the complexes are acyclic.

The cohomology of the complexes  $\text{K}^{x, y}_{2k+1}$,
$\text{K}^{y, x}_{2k+1}$ is one-dimensional in case of even $k$, and
elements $\underline{h}_{2k+1}$, $\overline{h}_{2k+1}$ are a
representatives of the corresponding cohomology class, respectively.
If $k$ is odd, then the complexes are acyclic.
\end{lemma}

We prove Lemma 3 and Proposition 4  for the odd case by induction.

The reader is encouraged to check the induction base for $K^{\alpha,
\beta}_{2k+1}$, $k\leq3$, $K_{m}$, $m\leq3$.

Applying simultaneous induction, we use the following
considerations.

There exist an exact triple (because of the monomial basis, see
Lemma 2):
$$ 0 \longrightarrow  \bigoplus\limits_{\alpha\neq\square
\text{ and } \beta\neq\square } K^{\alpha, \beta}_{2k+3}
\longrightarrow  K_{2k+3} \longrightarrow
\bigoplus\limits_{\alpha=\square \text{ or } \beta=\square }
K^{\alpha, \beta}_{2k+3}  \longrightarrow  0$$

As is well known, a long exact sequence of cohomologies is
associated with exact triple. But due to the inductive assumption,
we obtain the next exact sequence for $K_{2k+3}$:
$$0 \longrightarrow H^{k+1}(K_{2k+3}) \longrightarrow H^{k+1}(K_{2k+1}) \longrightarrow H^{k+2}(\bigoplus\limits_{\alpha\neq\square \text{ and } \beta\neq\square }
K^{\alpha, \beta}_{2k+3}) \longrightarrow
H^{k+2}(K_{2k+3})\longrightarrow 0$$

Similar exact sequences can also be written for all $K^{\alpha,
\beta}_{2k+3}.$

Because we know Euler characteristics of all $K^{\alpha,
\beta}_{2l+1}$ and $K_{2l+1}$ itself for any $l$ (absolute values do
not exceed 1) from paper $[1]$ and thanks to Lemma 2, then, to prove
the simultaneous induction, we just need to show that
$$(x_0+y_0)\circ\underline{h}_{2k+1} \neq 0$$ in corresponding
``2k+3'' spaces.

Up to a change of basis, we have:

\begin{minipage}[h]{.45\linewidth}
\begin{equation*}
(x_0+y_0)\underline{h}_{2k+1}=
 \begin{cases}
     (x_0+y_0)x_{-k} \underline{h}_{2k-3} y_{k} = x_{-k} \underline{h}_{2k-1} y_{k}, &\text{$k \mod 2 \neq 0$,}\\
   (x_0+y_0)y_{-k} \underline{h}_{2k-3} y_{k} = y_{-k} \underline{h}_{2k-1} y_{k}, &\text{$k \mod 2 = 0$.}
 \end{cases}
\end{equation*}
\end{minipage}

An element $x_{-k} \underline{h}_{2k-1}
y_{k}\;(y_{-k}\underline{h}_{2k-1} y_{k})$ contains only one
``prohibited'' pair:
$$x_{-k}y_{-k+1}\;(y_{-k}x_{-k+1})$$

It follows from the relations $\verb"REL"$ in $K_{2k+3}$ that
$$x_{-k}y_{-k+1} = -x_{-k+2}y_{-k-1},\;$$
$$y_{-k}x_{-k+1}=-y_{-k+2}x_{-k-1}.$$

Using the latter equations, it is possible to rewrite  $x_{-k}
\underline{h}_{2k-1} y_{k}\;(y_{-k}\underline{h}_{2k-1} y_{k})$ as
an element of the monomial basis, therefore, it is not equal to
zero, see Lemma 2.

The statement for  $K_{2k}$ follows from the statement for
$K_{2k-1}$ and Lemma 3, as it is possible to write exact triple for
$K_{2k}$ with $K_{2k-1}$ and $K^{x, \square}_{2k-1}, K^{y,
\square}_{2k-1}._\bigtriangledown$

\begin{example} $(x_{0}+y_{0})(x_{-1}y_{1}+x_{1}y_{-1})\neq0  \text{ in } K_{5}$

Let's use the relation
$x_{-2}y_{2}+x_{2}y_{-2}+x_{0}y_{0}=-(x_{-1}y_{1}+x_{1}y_{-1}):$
$$-(x_{0}+y_{0})(x_{-1}y_{1}+x_{1}y_{-1})=(x_{0}+y_{0})(x_{-2}y_{2}+x_{2}y_{-2})\neq0.$$
\end{example}

\section{Appendix: a generalization}

Conjecture is that it is possible to universalize the ideas for some
m-leg ladder models. That way, we obtain representations of $N=2,
\dots$ algebras.

\subsection{Graded Euler characteristic}

Let's recall basic definitions.

For each graph $\Gamma$ we can construct a statistical model in
which the set of configurations is the set of arrangements of
particles at graph vertices such that at each vertex at most one
particle is located and two particles cannot be located at vertices
joined by an edge.

To any graph $\Gamma$ there corresponds the fermion algebra
$A(\Gamma)$ defined as follows.

Let $S$ be the vertex set of $\Gamma$, and let $T \subset S\times S$
be its edge set. The algebra $A(\Gamma)$ is generated by $\psi_s, s
\in S$. The defining relations are
$$\psi_{s_1}\psi_{s_2}+\psi_{s_2}\psi_{s_1}=0,\;
\psi_{s_1}\psi_{s_2}=0, (s_1, s_2) \in T.$$

From $A(\Gamma)$ we can construct a complex $K(\Gamma)\cong
A(\Gamma)$ with differential $\sum\limits_{s \in S}\psi_{s}$.

The statistical sum of the model is the sum over the space of all
configurations. The contribution of each configuration depends on
parameters. For some parameters values, the contribution of any
configuration to the statistical sum equals $\pm 1$, depending on
the parity of the number of particles in the given configuration.
This statistical sum is naturally interpreted as the Euler
characteristic of the complex $K(\Gamma)$.

By a weight system $\{w\}$ on the graph $\Gamma$ we mean a function
on $S$ assigning a number $w(s)$ to each vertex $s$. For each
configuration $\chi = \{s_1, s_2, \dots, s_l\}$ we define its
"energy"
$$E^w(\chi)= (-1)^l q^{\sum\limits_{i=1}^l w(s_i)}, q \in
\mathbb{C}^{*}.$$

We introduce the statistical sum $\sum E^w(\chi)$, where the
summation is over all configurations. We refer to this as the graded
Euler characteristic and denote it by
$\text{E}_q^{\text{w}}(\Gamma)$.

\subsection{Cyclic 3-leg triangular
ladder model}

Cyclic 3-leg triangular ladder is the 4-leg triangular ladder, where
first and last rows coincide (are glued), see Fig. 5.

\begin{figure}[h]
\centering
\includegraphics[width=10.0cm]{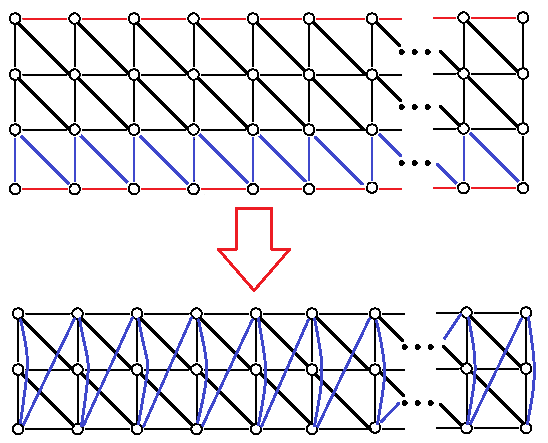}
\caption{cyclic 3-leg triangular ladder with $n$ columns}
\label{fig.0}
\end{figure}

We denote the graph from Fig. 5 by $\Gamma^3_n$.

Let $E^3_n$ denote the Euler characteristic of the complex
$K(\Gamma^3_n)$.

\begin{proposition} $E^3_{n} = (E^3_{1})^{[\frac{n+1}{2}]}= (-2)^{[\frac{n+1}{2}]}$, where $[ * ]$ denotes the integer part
of a number.\end{proposition} \textit{Proof. }It is easy to see that
such relation holds:
$$E^3_n = E^3_1 E^3_{n-2},$$ as any non-empty arrangement of
particles in second column of Fig. 5 partitions the graph into two
disconnected graphs, one of which is a one-point space. The Euler
characteristic of a disconnected graph is the product of the Euler
characteristics of its connected components.

In addition, $E^3_1=-2._\bigtriangledown$

Let's consider the weight system (which we denote by $f$), where
 each fermion in the column number $i$, $i=1, \dots, n,$ has weight
$[\frac{n}{2}]+i-1$.

To the cyclic 3-leg triangular ladder with $n$ columns we assign a
table of weights and numbers of fermions; namely, each cell of this
table contains the number of all admissible arrangements of a given
number of fermions with given weights.

The number of a column in the table corresponds to the number of
arranged fermions. The first column corresponds to $0$ fermions; it
is impossible to arrange more than $n$ fermions. The rows correspond
to weights. We assume that the configuration with no fermions has
weight $0$.

Below we give examples of such tables. We leave a cell empty if
there exist no configurations with given weight and given number of
fermions.

\begin{figure}[h]
\centering
\includegraphics[width=16.0cm]{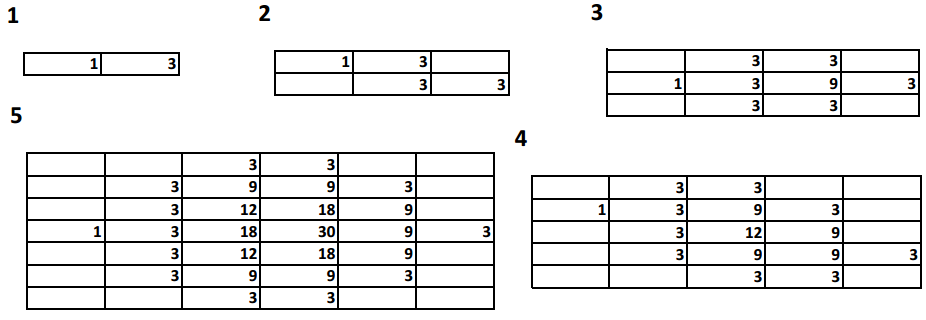}
\caption{first tables} \label{fig.0}
\end{figure}

For each $j =0, \dots ,n,$ to the $j$th column of a table there
naturally corresponds a Laurent polynomial in the independent
variable $q$. We denote these polynomials by $g^j_n(q)$. For
example, $g^3_3(q)= 3q+9+3q^{-1}$.

Clearly, $\text{E}^{f}_q(\Gamma^3_n) = \sum\limits_{i=0}^n
(-1)^{i}g_n^j(q).$

\begin{proposition} The relation $\text{E}^{f}_q(\Gamma^3_n)  = E^3_{n}$ holds.\end{proposition}
\textit{Proof.} This can be proved in the same manner as in paper
$[1]._\bigtriangledown$

All definition connected with semi-infinite configurations are the
same. The reader is encouraged to check it.

\begin{proposition} $$\widehat{F}(t, q)= \prod\limits_{m=1}^{\infty} \frac{1+2q^m}{1-q^m}\sum\limits_{n=-\infty}^{\infty}
 q^{n(n-1)/2} t^{n},$$ where $\widehat{F}(t, q)$ means the statistical
sum of the set of semi-infinite configurations in case of the cyclic
3-leg triangular ladder with the fixed type of tail relative to $C$
and $U._\bigtriangledown$ \end{proposition}

The fermion algebra for the infinite in both directions graph as per
Fig. 5 is the following algebra of anti-commuting elements $x_i$,
$y_i$ and $z_i$, $i \in \mathbb{Z}$:
$$\mathbb{C}[\dots, x_{-1}, y_{-1}, z_{-1}, x_{0}, y_{0}, z_{0},  x_{1}, y_{1}, z_{1}, \dots]$$
with relations: $$x_iy_i=0,\; x_iy_{i+1}=0;$$
$$y_iz_i=0,\; y_iz_{i+1}=0;$$
$$z_ix_i=0,\; z_ix_{i+1}=0;$$
$$x_ix_{i+1}=0,\; y_iy_{i+1}=0,\; z_iz_{i+1}=0.$$

Let's denote generating functions $$X(t) = \sum_{i\in  \mathbb{Z}}
x_i t^{-i},\;Y(t) = \sum_{i\in \mathbb{Z}} y_i t^{-i},\;Z(t) =
\sum_{i\in \mathbb{Z}} Z_i t^{-i}.$$

Let's determine the deformation of the fermion algebra for the
infinite in both directions graph as per Fig. 5. It is the algebra,
generated by anti-commuting elements $x_i$, $y_i$ and $z_i$,
satisfying relations below:
$$X(t)Y(t) = 0,\;Y(t)Z(t) = 0,\;Z(t)X(t) = 0,$$
$$X(t)X(-t)= 0,\;Y(t)Y(-t) = 0,\;Z(t)Z(-t) = 0.$$

The deformation preserves dimensions in the corresponding induced
representations with extreme vectors.

\subsection{Another model}

Let's consider such graph (which we denote by
$\widetilde{\Gamma}^m_n$):
\begin{figure}[h]
\centering
\includegraphics[width=10.5cm]{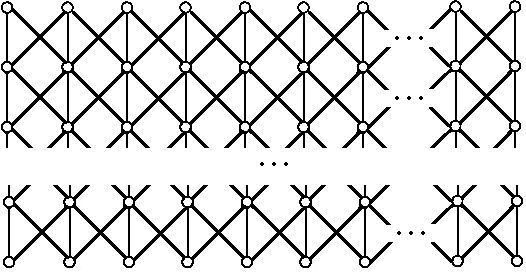}
\caption{m-leg ladder with $n$ columns} \label{fig.0}
\end{figure}

Let $\widetilde{E}^m_n$ denote the Euler characteristic of the
complex $K(\widetilde{\Gamma}^m_n)$.

\begin{proposition}$\widetilde{E}^m_{n} = (\widetilde{E}^m_{1})^{[\frac{n+1}{2}]}$, where $[ * ]$ denotes the integer part
of a number.\end{proposition} \textit{Proof. }It is easy to see that
such relation holds: $$\widetilde{E}^m_n = \widetilde{E}^m_1
\widetilde{E}^m_{n-2},$$ by the same reasons as in Proposition 5.

Moreover,
$$\widetilde{E}^m_{1}=-\widetilde{E}^{m-3}_{1},$$ where
$$\widetilde{E}^1_{1}=0,\;\widetilde{E}^2_{1}=-1,\;\widetilde{E}^3_{1}=-1._\bigtriangledown$$

Let's recall the weight system (which is denoted by $f$), where
 each fermion in the column number $i$, $i=1, \dots, n,$ has weight
$[\frac{n}{2}]+i-1$.

\begin{proposition}$E^{f}_q(\widetilde{\Gamma}^3_1)  = \widetilde{E}^3_{1},$ $E^{f}_q(\widetilde{\Gamma}^3_2) = \widetilde{E}^3_{2},$ $E^{f}_q(\widetilde{\Gamma}^3_3) = 3-q-q^{-1}  \neq \widetilde{E}^3_{3}.$\end{proposition}
\textit{Proof.} A direct calculation$._\bigtriangledown$

\begin{remark}Proposition 9 does not mean that the idea does not work due to Proposition 5.\end{remark}

Let $\Theta$ denote such set of weight systems, for which is assumed
that each fermion in the column number $1$ on Fig. 7 has weight $0$.

\begin{proposition} The relation $E^{\text{o}}_q(\widetilde{\Gamma}^{3m+1}_n)
= \widetilde{E}^{3m+1}_n$ holds for any $o \in \Theta.$
\end{proposition} \textit{Proof. }It follows from Proposition 8, as
$\widetilde{E}^{3m+1}_1 = 0._\bigtriangledown$

Higher School of Economics -- National Research University, Moscow\\
\small\texttt{VvS@myself.com}
\end{document}